\documentclass[graybox,12pt]{svmult}

\usepackage{mathptmx}       
\usepackage{helvet}         
\usepackage{courier}        

\usepackage{makeidx}         
\usepackage{graphicx}        
\usepackage{multicol}        

\usepackage{amsmath}
\usepackage{amsfonts}
\usepackage{color}     
\usepackage{url}
\usepackage{verbatim}
\usepackage{enumerate}

\makeindex             

\usepackage{amssymb}
\usepackage{amsmath}
\usepackage{epsfig}
\usepackage{mathptmx}       
\usepackage{helvet}         
\usepackage{courier}        
\usepackage{type1cm}        

\usepackage{makeidx}         
\usepackage{graphicx}        
\usepackage{multicol}        
\usepackage[bottom]{footmisc}

\usepackage{amsbsy}
\usepackage{MnSymbol}
\usepackage{amsfonts}
\usepackage{amsopn}
\usepackage{dsfont}
\usepackage[english]{babel}
\usepackage[applemac]{inputenc}
\usepackage[colorlinks=true,breaklinks=true]{hyperref}
\usepackage[T1]{fontenc}
\usepackage{xcolor}
\usepackage{enumerate}
\usepackage{latexsym,amsmath,eufrak}
\usepackage{xcolor}
\usepackage{stackrel}
\newenvironment{proofp}{\noindent {\it Proof. }}{\ \ \ $\square$\\}
\newenvironment{proofa}{\noindent {\bf Proof.}}{\ \ \ $\Box$}

\newcommand{\limas}{\stackrel{a.s.}{\longrightarrow}}
\newcommand{\limdis}{\stackrel{d}{\longrightarrow}}

\def\indi{\mathrm{\hspace{0.2em}l\hspace{-0.55em}1}}

\newcommand{\Exp}{\mathbb{E}}
\newcommand{\Vari}{\mathbb{V}}
\newcommand{\Prob}{\mathbb{P}}
\newcommand{\eqdis}{\stackrel{d}{=}}
\newtheorem{teo}{Theorem}[section]
\newtheorem{prop}[teo]{Proposition}
\newtheorem{lem}[teo]{Lemma}
\newtheorem{exaample}[teo]{Example}
\newtheorem{rmq}[teo]{Remark}
\newtheorem{cor}[teo]{Corollary}

\begin{document}
\title*{Unbalanced multi-drawing urn with random addition matrix II}
\author{
Rafik Aguech
\thanks{Department of mathematics , University of Monastir, Tunisia, \email{rafik.aguech@ipeit.rnu.tn},}  \thanks{University of Monastir, Tunisia, \email{selmiolfa3@yahoo.fr}},
Wissem Jedidi$^\star$
\thanks{Universit\'e de Tunis El Manar, Facult\'e des Sciences de Tunis, D\'epartement de Math\'ematiques, Laboratoire d'Analyse Math\'ematiques et Applications LR11ES11.  2092 - El Manar I, Tunis, Tunisia \email{ wissem.jedidi@fst.utm.tn}},
Olfa Selmi$^{\dagger}$
}
\maketitle
\abstract{At at stage $n \geq 1$, we pick out  at random  $m$ balls, say $k$ white balls and $m-k$ black balls. We inspect the colors and then we return the balls, according to a predefined  replacement matrix, together with $(m-k)\; X_n$ white balls and  $k\; Y_n$ black balls.  Under the  assumption that the seqence $(X_n, \;Y_n)$ is bounded and i.i.d.   Aguech \& Selmi \cite{R&O} proved a strong law of large numbers (SLLN) and a central limit theorem (CLT) on the proportion of white balls. Here, we extend the last results as follows.  In a first step, we obtain the same results  under the assumption  that $(X_n, \;Y_n)$ have the same distribution and that $\Exp[X_1^2+Y_1^2]<\infty$.  In a second step, removing the assumption of identical distribution  and  assuming that   $\Exp[X_n^2+Y_n^2] <\infty$ for all $n$,  we prove the  SLLN and CLT for the total number of balls in the urn.}


\section{Introduction}
The following notations will be used throughout the paper. The indicator function
is denoted by $ \mathbf {1}_{\{.\}}$, the abbreviation $a.s.$ stands for almost surely and the notations $\limas$, $\limdis$ respectively stand for  the almost surely convergence and the convergence in distribution. For a sequence of i.i.d random variables ${(R_n)_n}$, we denote by
\[  \mu_R=\Exp[{R_1}]\quad  \mbox{ and }\quad \sigma_R^2=\Vari[R_1], \]
respectively the mean and the variance of a random variable $R_1$. For two sequences of real number $x_n$ and $y_n$, we use  the classical Landau notations $x_n\sim y_n$ if $\lim_n x_n/y_n=1$, $x_n=o(y_n)$ if $\lim_n x_n/y_n=0$ and $x_n=O(y_n)$ if $\limsup_n |x_n/y_n|<\infty$. The mention $a.s.$ is added if the sequences are random.\\

The P\'olya urn is a simple and powerful model that still widely used since its appearance in  P\'olya and Eggenberger \cite{Polya}.   Historically, the classical model is an urn containing white balls and blue balls. At each discrete time $n\geq 1$, we pick out, uniformly at random,   a ball, after inspecting its color, the composition of the urn is then evolved according to a replacement rule represented by a matrix given by
\begin{center}
$R=\left(
\begin{array}{cc}
                  a & b \\
                  c & d \\
                \end{array}
              \right),$
\end{center}
where $a,b,c$ and $d$ are integers. Since the mechanism of such a Markov chain is simple, we find its applications in several fields such as computer sciences, diseases spreading, finance, random research trees... For more details about the history and the diversity of the applications of the urns processes we refer the reader to \cite{Kotz&jon,  M.Hist}\\

\ Afterwards, a generalized urn model has received interest in the literature, this process evolves in the following way: we fix an integer $m\geq 1$, at each discrete time step we draw randomly $m$ balls from the urn (either with or without replacement), the colors are inspected and the balls are reinserted to the urn
together with other balls, depending on the sampled balls and the addition rule is a matrix $R$ given by \[R=\left(%
\begin{array}{cc}
  a_0 & b_0 \\
  \vdots& \vdots \\
  a_{m-1} & b_{m-1} \\
  a_m & b_{m-1} \\
\end{array}%
\right).\]

Such a model was firstly introduced and studied by Chen \& Wei \cite{chen&wei}, Chen \& Kuba
\cite{chen&kuba} and  Kuba, Mahmoud \& Panholzer \cite{kuba&mah&al}, who considered a balanced model  i.e.
$$a_i+b_i>0, \;\; \mbox{is a constant for all}\;  0\leq i\leq m$$
for both cases
$$a_i=c(m-i)  \quad \mbox{and}\quad a_i=c\,i,\quad \;1\leq i\leq m , \quad \mbox{(for some integer}\; c).$$
This model was further developed by  Kuba \&  Mahmoud \cite{Kuba&mahmoud, Kuba-Mahmoud2} and  Kuba \&  Sulzbach \cite{K&Sulz}, who studied a more general model under the assumptions of balance (i.e. $a_i+b_i=\sigma>0$) and affinity (the conditional expectation of the number of white balls in the urn satisfies an affine relation). They proved  that the study can be reduced to the case when the coefficient of the replacement matrix satisfies itself the affinity condition given by
$$a_i=a_0+\frac{a_m-a_0}{m} i, \quad \mbox{for}\; 1\leq i\leq m-1.$$

\ The works of   Aguech,  Lasmar \&  Selmi \cite{A&L&S} and  Lasmar, Mailler \&  Selmi \cite{L&M&S}  focused on this model by removing some assumptions. The challenge of the pre-cited works, was to give a description of the behavior of the urn's composition, while the urn is neither  balanced nor affine. Some of their  adapted stochastic approximation methods also gave consistent results.\\

Multiple drawing urns with random addition matrix have also been studied by Aguech \&  Selmi \cite{R&O}, who considered the model defined as follows:
\begin{itemize}
\item We start with an urn containing white and black balls. Then, we fix an integer $m\geq 1$ and   suppose that at the beginning the urn contains more than $m$ balls (to ensure that the first draw is possible).

\item At each discrete time step $n\geq 1$, we pick out  at random $m$ balls from the urn. After inspecting the color of the sample, the balls are replaced to the urn.

\item  If the sample is composed by $k$ black balls and $m-k$ white balls, then
\begin{equation}\label{xy}
\mbox{\it we add  $(m-k)\;X_n$ white balls and $k\;Y_n$ black balls, where $X_n, \; Y_n=0,1,2, \ldots$}
\end{equation}
The sequences $(X_n)_n$ and $(Y_n)_n$ are independent. If $R_n=X_n$ (respectively  $R_n=Y_n$), then
    \begin{equation} \label{ASG}
    \mbox{\it the r.v.'s $R_n$  are independent copies of some integer-valued random variables $R=X$ (respectively  $R=Y$).}
    \end{equation}
\item We start with two given integers $W_0,\;B_0$. In order to ensure that the first draw  is possible, the initial composition  of the urn should satisfy
    \begin{equation} \label{t0}
     T_0:=W_0+B_0\geq m
     \end{equation}
     The composition of the urn, at a time $n\geq 1$, is then given by a vector $(W_n,B_n)$, where
     \begin{eqnarray}
     T_n&:=& \mbox{\it the total number of balls in the urn}, \label{tn}\\
     W_n&:=& \mbox{\it the number of white balls} \in \{0,1,\ldots, T_n\}, \quad B_n :=  \mbox{\it the number of black balls,}=T_n-W_n. \label{wbn}
         \end{eqnarray}
     It follows that, conditionally on $(W_{n-1},B_{n-1)}$, the probability of drawing $k$ white balls and $(m-k)$ black balls,  is for  given, each time $n\geq 1$, by
     $$\frac{ {{W_{n-1}}\choose{k}}\; {{B_{n-1}}\choose{m-k}} }{ {{T_{n-1}}\choose{m}} },\quad \mbox{\it where}\; {{j}\choose{i}}= \frac{j!}{i!\,(j-i)!}.$$
\item We denote by
\begin{eqnarray}
\xi_n&:=&\mbox{\it the number of white ball among the $n^{th}$ sample $\in\{0,1,\ldots m\}$,} \label{xn}\\
Z_n&:=&\frac{W_n}{T_n} =\mbox{\it the proportion of white balls in the urn after $n$ draws $\in [0,1]$.} \label{zn}
\end{eqnarray}
Note that the r.v.'s $T_n$, $W_n$  and $Z_n$ are not deterministic, hence they are of primary  interest for studying.\\

Let $\mathcal{F}_n$ be the $\sigma$-field generated by the first $n$ draws.  Then, conditionally on $\mathcal{F}_{n+1}$, the random variable $\xi_n$ in \eqref{xn} follows the hypergeometric distribution with parameters ($m,W_{n-1},T_{n-1}$). Thus, the mean and the variance of $\xi_n$ are given by
\begin{equation}\label{hypo}
\Exp[\xi_n|\mathcal{F}_{n-1}]=m\, \frac{W_{n-1}}{T_{n-1}}=m \; Z_{n-1},\quad \Vari[\xi_n|\mathcal{F}_{n-1}]=m\frac{W_{n-1}}{T_{n-1}}\left(1-\frac{T_{n-1}-m}{T_{n-1}-1}\right).
\end{equation}
Conditioning on the update of the urn at the time $n$, the number of balls in the urn satisfies the recursive equation
\begin{equation}\label{recusrion-opp}\left(%
\begin{array}{c}
  W_{n+1} \\
  B_{n+1} \\
\end{array}%
\right)=\left(%
\begin{array}{c}
  W_n \\
  B_n \\
\end{array}%
\right)+\left(%
\begin{array}{cc}
  0 & X_{n+1} \\
  Y_{n+1} & 0 \\
\end{array}%
\right)\left(\begin{array}{c}
  m-\xi_{n+1}\\
  \xi_{n+1}\\
\end{array} \right),\quad n\geq 0,
\end{equation}
and in view of the evolution of the urn process,  the total number of balls in the urn after $n$ draws  satisfies
\begin{equation} \label{recursion-T_n}
T_n=T_0+\sum_{i=1}^nX_i(m-\xi_i)+\sum_{i=1}^n\xi_iY_i, \quad \mbox{ a.s}.
\end{equation}
\end{itemize}

\noindent Under  a condition  stronger than \eqref{ASG}, namely
\begin{equation} \label{ASc}
\mbox{\it  the sequence  $R_n=X_n$ (respectively  $R_n=Y_n$)   is i.i.d and $R_n$ is bounded,}
\end{equation}
Aguech \&  Selmi \cite{R&O} showed that  there is a  SLLN and a CLT for $Z_n$ and for $T_n/n$. The contribution of this paper is to extend the their model. For instance, in Section \ref{Model 1}, we remove  the stationary distribution and boundedness conditions \eqref{ASc} and   we obtain the same results if we assume for   $R_n=X_n$ (respectively  $R_n=Y_n$):
\begin{equation} \label{A}
\mbox{\it  the sequence  $R_n$  is i.i.d , $R_n$ is bounded below by some  positive constant $L$, and is square integrable}.
\end{equation}
In Section \ref{Model 2}, we extend the last model by removing the assumption of identical distribution. We assume instead that
\begin{equation} \label{B}
\mbox{\it  the sequence  $R_n$  is independent and $\Exp[X_n^2+Y_n^2] <\infty$, for all $n\geq 1$,}
\end{equation}
and we prove the  SLLN and CLT with another normalization of the total number of balls  in the urn $T_n$.
We also obtain the almost sure asymptotic behavior of $T_n$ and  $Z_n$  under the additional assumption:
that the sequence $\Exp[X_n]$ and  $\Vari[X_n]$ are regularly of orders  $\alpha, \; \gamma>-1$ , i.e.
$$\mbox{\it  the sequence $\Exp[X_n]$ and  $\Vari[X_n]$ are regularly of orders  $\alpha, \; \gamma>-1$},$$
i.e.
\begin{equation}\label{D}
\Exp[X_n] = n^{\alpha}\, l_1(n), \quad  \Vari[X_n] = n^{\gamma} l_2(n), \quad l_i, \;i=1,2 \; \;\mbox{are a slowly varying functions}.
\end{equation}
(See \eqref{slow} below for slowly varying functions).\\

\noindent Section \ref{appendix} is an appendix, where we recall some useful results.
\section{Results under assumption \eqref{A}}\label{Model 1}
In this seection, we provide a SLLN and  CLT for the proportion  of white balls and of the total number of balls in the urn.
\subsection{The SLLN for $Z_n$ and $T_n$}\label{SLLN}

\begin{lem}\label{sum-T_n} Under assumption \eqref{A}, the total number $T_n$ defined in \eqref{tn}, satisfies
\[ \sum_{n\geq 0} \frac{1}{T_n^2}<\infty\quad \text{and}\quad \sum_{n\geq
0} \frac{1}{T_n}=\infty, \quad \mbox{a.s}.
\]
\end{lem}
\begin{proofp} From \eqref{recursion-T_n}, we have $T_n\geq T_0+m \; n\; \min_{1\leq i\leq n}(X_i,Y_i)  \geq T_0+ m\; n\; L$, then,
$$\sum_n\frac{1}{T_n^{2}}\leq \sum_n\frac{1}{(T_0+m\; n\; L)^2}< +\infty, \quad a.s.$$
On the other hand, since   $0\leq\xi_i\leq m$ for all $i\geq 1$, then
$$\frac{T_n}{n}\leq \frac{T_0}{n}+\frac{m}{n}\sum_{i=1}^{n}(X_i+Y_i), \quad a.s.$$
By the SLLN,  for all $\varepsilon >0$ and, a.s. for all $\omega \in \Omega$,  there exists $n_0(\omega)>0$ such that
$$\frac{T_n(\omega)}{n}\leq \varepsilon+ (\mu_X+\mu_Y)\, m, \quad n\geq n_0(\omega).$$
The latter entails that for all $n\geq n(\omega)$
$$\frac{1}{T_n(\omega)}\geq \frac{1}{n{\left(\varepsilon +\left(\mu_X+\mu_Y\right)m\right)}},$$ and consequently, $\sum_n 1/{T_n(\omega)} =\infty$.
\end{proofp}
\begin{prop}\label{prop-cv} Under the assumption  \eqref{A}, the proportion $Z_n$ of white balls after $n$ draws  satisfies
\begin{equation}\label{zz}
Z_n \limas z_\star:=\frac{\sqrt{\mu_X}}{\sqrt{\mu_X}+\sqrt{\mu_Y}},\quad a.s.
\end{equation}\end{prop}
\begin{proofp} Recall that, in view of the recursive equation (\ref{recusrion-opp}), the number of white balls $W_n$ and the total number of balls in the urn after $n$ draws satisfy
$$W_{n+1}=W_n+X_{n+1}(m-\xi_{n+1})\quad \text{and}\quad T_{n+1}=T_n+m\,X_{n+1}+\xi_{n+1}(Y_{n+1}-X_{n+1}).$$
Then, the proportion $Z_n$  can be described via the  stochastic algorithm
\begin{equation}\label{rec-Z_n}
Z_{n+1}=Z_n+\gamma_{n+1}\big(f(Z_n)+\epsilon_{n+1}\big),
\end{equation}
where
\begin{equation} \label{bruit}
f(z)=m\; \big((\mu_X-\mu_Y)\; z^2-2\;\mu_X \;z+\mu_X\big), \quad \epsilon_{n+1}:= D_{n+1}-\Exp[D_{n+1}|\mathcal{F}_n],
\end{equation}
and
$$D_{n+1}:=\xi_{n+1}\,\left[Z_n\;(X_{n+1}-Y_{n+1})-X_{n+1}\right]+m\,X_{n+1}\,(1-Z_n).$$
Consider the recursion (\ref{rec-Z_n}) in the form of \eqref{eq-monro} in the Appendix, with
\begin{equation} \label{agf}
\gamma_n:=\frac{1}{T_{n}}\quad \text{and}\quad f(z):=m \; \big((\mu_X-\mu_Y)\;z^2-2\;\mu_X \; z+\mu_X\big).
\end{equation}
Let us check the conditions of Theorem \ref{Monro}. Recall that in \eqref{t0}, we have $\gamma_0\leq 1/m$. We also have:
\begin{enumerate}[{\it (a)}]
\item The first condition of Theorem \ref{Monro} is satisfied since
      he r.v.'s $Z_n,\; \epsilon_n$ are bounded, $\Exp[\epsilon_{n+1}|\mathcal{F}_n]=0$, $\gamma_n$ is nonnegative, decreases towards zero, and $\gamma_0=1/{T_0}>0$ is a deterministic constant.
\item  Since $f$ is a second order polynomial, we easily  verify that the constant $z_\star$, given by \eqref{zz}, is the unique solution of $f(z)=0$ and $f'(z)<0$.
\item The constant $K$ of the condition {\it (c)} in Theorem \ref{Monro} corresponds to $K := m \; \Big(|\mu_X-\mu_Y|+2\;\mu_X \Big)$, since for all $z\in [0,1]$,
     $$\frac{|f(z)|}{m}\leq |\mu_X-\mu_Y|\;z^2+2\;\mu_X\; z+\mu_X \leq \Big(|\mu_X-\mu_Y|+2\;\mu_X \Big)\;z+ \mu_X \leq  \Big(|\mu_X-\mu_Y|+2\;\mu_X\Big)(1+z).$$
\item Using the expression of $\epsilon_{n+1}$ in \eqref{bruit}, the facts  $X_{n+1}$ and $Y_{n+1}$ are independent of $\mathcal{F}_n$, we observe that
     $$\epsilon_{n+1}=\xi_{n+1}\Big(Z_n\,(X_{n+1}-Y_{n+1})- X_{n+1}\Big)+m\,X_{n+1}\,(1-Z_n)-m\,Z_n \,\Big(Z_n\,(\mu_X-\mu_Y)-\mu_X\Big)-m\mu_X\,(1-Z_n\,).$$
     By \eqref{hypo}, recall that $\Exp[\xi_{n+1}|\mathcal{F}_n]=m\,Z_n$, then
     \begin{eqnarray}
     \Exp\left[\epsilon_{n+1} ^2|\mathcal{F}_n\right]&=&\Exp\left[\xi_{n+1} ^2|\mathcal{F}_n\right]\; \Big(Z_n^2\left(\mu_{X^2}+\mu_{Y^2}-2\mu_X\mu_Y\right)+
     2\,Z_n\left(\mu_X\mu_Y-\mu_{X^2}\right)+\mu_{X^2}\Big)\nonumber  \\
     &&-m^2 Z_n^2\Big(Z_n\,(\mu_X-\mu_Y)-\mu_X\Big)^2 + m^2\,(1-Z_n)^2 \,\sigma_X^2
\label{rest} \end{eqnarray}
Finally, using the facts that $0\leq\xi_{n+1}\leq m$, $Z_n \in [0,1]$, and after elementary computations, we see that there exists a constant $C_{\epsilon}>0$ such that

\begin{equation} \label{dddd}
\Exp\left[\epsilon_{n+1} ^2|\mathcal{F}_n\right] \leq  C_{\epsilon} \Longrightarrow  \Exp[|\epsilon_{n+1}| \big|\mathcal{F}_n] \leq  \sqrt{C_{\epsilon}}.
\end{equation}
Lemma \ref{sum-T_n} and \eqref{dddd} ensure  that, $\gamma_n$ given in \eqref{agf}, satisfies
\[\sum_{n\geq 0}\gamma_n=\infty \quad \text{and}\quad \sum_{n\geq 0}\gamma_n^2\;
\Exp[\epsilon_{n+1}^2|\mathcal{F}_n] \leq C_{\epsilon}\sum_{n\geq 1}\gamma_n^2<\infty, \quad a.s.\]
\end{enumerate}
In conclusion, Theorem \ref{Monro} asserts  that $Z_n \limas z_\star$.
\end{proofp}

The SLLN on the proportion of  black balls of the urn is obtained by considering $1-Z_n$.
\begin{prop}\label{total-cv}  Under  assumption  \eqref{A}, the total number $T_n$ of balls and the number of white balls $W_n$ after $n$ draws, satisfy
\begin{equation} \frac{T_n}{n} \limas m\sqrt{\mu_X\mu_Y}\quad \mbox{and} \quad \frac{W_n}{n}  \limas  \frac{m\;\mu_X \;\sqrt{\mu_Y}}{\sqrt{\mu_X}+\sqrt{\mu_Y}}.
\end{equation}
\end{prop}
\begin{proofp} If $G_n:=\sum_{i=1}^n\Big(\xi_i(Y_i-X_i)-\Exp[\xi_i(Y_i-X_i)|\mathcal{F}_{i-1}]\Big),$ then  (\ref{recursion-T_n}) gives the representation
\begin{equation}
\label{eq-G_n}\frac{T_n}{n}=\frac{T_0}{n}+\frac{m}{n}\sum_{i=1}^n X_i +\frac{m(\mu_Y-\mu_X)}{n} \sum_{i=1}^nZ_{i-1}+\frac{G_n}{n}.
\end{equation}
Since $(X_i)_{i\geq 1}$ are i.i.d. random variables, the SLLN insures that $\lim_{n\rightarrow +\infty}\sum_{i=1}^n X_{i}/n= \mu_X$ and Ces\'aro's lemma yield
$$\sum_{i=1}^n\frac{Z_{i}}{n}  \limas z_\star.$$
It remains to  prove  $G_n/n \limas 0$. We proceed as follows: the sequence $(G_n,\mathcal{F}_n)_{n\geq 0}$ is a martingale difference  whose quadratic variation is given by
\begin{equation*}<G>_n:=\sum_{i=1}^n\Exp[\Delta G_i ^2|\mathcal{F}_{i-1}],\quad \mbox{where}\quad
\Delta G_n=G_n-G_{n-1}=\xi_n(Y_n-X_{n})-\Exp[\xi_n(Y_n-X_{n})\;|\;\mathcal{F}_{n-1}].\end{equation*}
A simple computation gives
$$\lim_{n\rightarrow +\infty}\Exp[\Delta G_n^2|\mathcal{F}_{n-1}]=\left(m\; z_\star\; (1-z_\star)+m^2\; {z_\star}^2\right)(\sigma_Y^2+\sigma_X^2),$$
where $z_\star$ is given in \eqref{zz}. Ces\'aro's lemma allows to conclude that
$$\lim_{n\rightarrow +\infty}\frac{<G>_n}{n} \limas \left(m\;z_\star \;(1-z_\star)+m^2\;{z_\star}^2\right) \Longrightarrow\frac{ G_n}{n}  \stackrel{a.s}{\longrightarrow} 0.$$
The last claim is immediate by writing $W_n/n= (W_n/T_n)(T_n/n)$.
\end{proofp}
\subsection{CLT for $Z_n$ and $W_n- T_n z_\star$}\label{CLT}
We give a CLT for the proportion of white balls in the urn.
\begin{teo}\label{last} Under the assumptions of Lemma \ref{sum-T_n}, the proportion of white balls in the urn after $n$ draws, $Z_n$, satisfies as $n$ tends  to infinity
\begin{equation} \label{tlc}
\sqrt{n}(Z_n-z_\star)\limdis\mathcal{N}\Big(0,\frac{P(z_\star)}{3\mu_X\;\mu_Y}\Big), \quad \mbox{where}\;\;P(z_\star)=\left(\sigma_X^2+\sigma_Y^2\right)\,z_\star^4-2\sigma_X^2\,z_\star^3+2\sigma_X^2\, z_\star^2-2\sigma_X^2\,z_\star+\sigma_X^2.
\end{equation}
The number, $W_n$, of white balls in the urn after $n$ draws satisfies
\begin{equation}
\sqrt{n}\left(\frac{W_n}{n}- \frac{T_n}{n}\;z_\star\right)\limdis\mathcal{N} \Big(0,\frac{m^2\, P(z_\star)}{3}\Big).
\end{equation}
\end{teo}
\begin{proofp} We recall that the sequences  $\epsilon_n, \;\gamma_n$ and the constant $C_\epsilon$  are respectively given in \eqref{bruit}, \eqref{agf}  and  \eqref{dddd}. To prove the claimed statement, we need to check the assumptions $(a),\ldots,(e)$ of  Theorem \ref{th1-renlund} in the Appendix, applied on the sequence $K_n:=\sqrt{n}(Z_n-z_\star)$. Observe that the representation  \eqref{rec-Z_n} could be rewritten in the form
\begin{equation}\label{rec-Z_n-zstar}
Z_{n+1}-z_\star=Z_n-z_\star+ \frac{1}{T_{n+1}}\big(f(Z_n)-f(z_\star)+\epsilon_{n+1}\big),
\end{equation}
and let us denote
\begin{equation}\label{gg}
K_n:=\sqrt{n}\left(Z_n-z_\star\right),\quad  \alpha_n=  \sqrt{1+\frac{1}{n}} -1 \quad \mbox{and} \quad
g(z):=\frac{f(z)}{z_\star-z}.
\end{equation}
Using the identity $f(Z_n)-f(z_\star) =g(Z_n)\left(z_\star-Z_n\right)$,   we rewrite \eqref{rec-Z_n-zstar} into a recursive equation on $K_n$  of the form \eqref{eq-renlund}:
$$K_{n+1}=\left(1-\frac{\Gamma_{n+1}}{n}\right)\,K_n+\frac{V_{n+1}}{\sqrt{n}},$$
where
\begin{equation}\label{gaga}
\Gamma_{n+1}:= \frac{n\;g(Z_n)}{T_{n+1}}- n\alpha_n \left(1-\frac{g(Z_n)}{T_{n+1}} \right)\quad \mbox{and} \quad  V_{n+1}:=\frac{\epsilon_{n+1}\sqrt{n\left(n+1\right)}}{T_{n+1}}.
\end{equation}

\begin{enumerate}[{\it (a)}]
\item Using \eqref{recursion-T_n}, we may decompose the r.v. $V_{n+1}$ in \eqref{gaga} in the form $$\sqrt{n}V_{n+1}=\frac{n\sqrt{n+1}}{T_{n+1}}\epsilon_{n+1} =n\;\sqrt{n+1}\;\epsilon_{n+1} \; \left(\frac{1}{T_{n+1}}-\frac{1}{T_{n}} +\frac{1}{T_n}\right)=n\;\sqrt{n+1}\;\epsilon_{n+1} \; \left(\frac{1}{T_n} - \frac{X_{n+1}\;(m-\xi_{n+1})+ Y_{n+1}\; \xi_{n+1}}{T_{n+1}\; T_{n}}\right),$$
    and obtain the inequalities
    $$-n\sqrt{n+1}\,m\;\left(X_{n+1}+Y_{n+1}\right)\;\frac{|\epsilon_{n+1}|}{T_{n+1}T_n}+\frac{n\sqrt{n+1}}{T_n}\epsilon_{n+1}\leq \sqrt{n}V_{n+1}\leq  n\sqrt{n+1}\, m\; \left(X_{n+1}+Y_{n+1}\right) \;\frac{|\epsilon_{n+1}|}{T_{n+1}T_n} +\frac{n\sqrt{n+1}}{T_n}\epsilon_{n+1}.$$
     Then,
    $$- \frac{\sqrt{n+1}\; m}{m^2\; n\; L^2} \; \left(X_{n+1}+Y_{n+1}\right) |\epsilon_{n+1}| +\frac{n\sqrt{n+1}}{T_n}\epsilon_{n+1}\leq \sqrt{n}V_{n+1}\leq \frac{\sqrt{n+1}\; m}{m^2\; n\; L^2}\;
    \left(X_{n+1}+Y_{n+1}\right) |\epsilon_{n+1}|+\frac{n\sqrt{n+1}}{T_n}\epsilon_{n+1}.$$
    In one hand, since $X_{n+1}$ and $Y_{n+1}$ are independent of $\epsilon_{n+1}$ and $\mathcal{F}_n$ and using \eqref{dddd}, we obtain
    $$\Exp\left[\frac{\sqrt{n+1}\; m}{m^2\; n\; L^2} \left(X_{n+1}+Y_{n+1}\right) |\epsilon_{n+1}| \;\big|\;\mathcal{F}_n\right]= \frac{\sqrt{n+1}\; m\big(\mu_X+\mu_Y\big)}{m^2\; n\;L^2}\Exp\left[|\epsilon_{n+1}| \; \big|\;\mathcal{F}_n\right]\limas 0.$$
    On the other hand,
    $$\Exp\left[\frac{n\sqrt{n+1}}{T_n}\epsilon_{n+1}\big|\mathcal{F}_n\right]= \frac{n\sqrt{n+1}}{T_n}\; \Exp\Big[\epsilon_{n+1} \big|\mathcal{F}_n\Big]=0, \quad  \mbox{for all} \; n\geq 1.$$
    All in one, we conclude that
    $$\Exp\Big[V_{n+1}|\mathcal{F}_n\Big]=o(n^{-1/2}), \quad a.s.$$
\item As in a), by \eqref{recursion-T_n}, we have  $T_{n+1} \geq m\;(n+1)\; L$, and then by \eqref{dddd},
      $$\Exp[V_{n+1}^2|\mathcal{F}_n]=\frac{n+1}{n} \Exp\left[\left(\frac{n+1}{T_{n+1}}\right)^2 \epsilon_{n+1}^2|\mathcal{F}_n\right] \leq \frac{3}{2\;m^2\; L^2}\;\Exp\left[\epsilon_{n+1}^2 \big|\mathcal{F}_n\right] \leq C:=\frac{3\; C_{\epsilon}}{2 \;m^2\; L^2} , \quad n\geq 1.$$
\item The sequence $(n+1)/T_{n+1}$ is bounded by $(mL)^{-1}$, and by Proposition \ref{total-cv}, it converges a.s. to $(m\,\sqrt{\mu_X\mu_Y})^{-1}$.  Then, using \eqref{rest},  we obtain
     $$ \Exp[V_{n+1}^2|\mathcal{F}_n]= \frac{n+1}{n}\,\Exp\left[\left(\frac{n+1}{T_{n+1}}\; \epsilon_{n+1}\right)^2|\mathcal{F}_n\right]\limas \frac{1}{m^2\mu_X\mu_Y}\lim_{n\rightarrow+\infty} \Exp\left[\epsilon_{n+1}^2|\mathcal{F}_n\right] =: \sigma^2>0,$$
     where
    \begin{eqnarray*}
    \sigma^2&=&\frac{z_\star^2}{\mu_X\mu_Y}\left[z_\star^2\left(\mu_{X^2}+
     \mu_{Y^2}-2\mu_X\mu_Y\right) +2z_\star\left(\mu_X\mu_Y-\mu_{X^2}\right)+\mu_{X^2}\right]\nonumber \\
      &&-\frac{z_\star^2}{\mu_X\mu_Y}\left[z_\star^2\left(\mu_X-\mu_Y\right)^2-2z_\star\mu_X\left(\mu_X-\mu_Y\right)+\mu_X^2 \right]+ \left(1-z_\star\right)^2\frac{\sigma_X^2}{\mu_X\mu_Y}.\label{sigma}
      \end{eqnarray*}
    Rearranging the last terms, we obtain the expression of $P\left(z_\star\right)$:
    \begin{equation}
     \mu_X\,\mu_Y\,\sigma^2=\left(\sigma_X^2+\sigma_Y^2\right)\,z_\star^4-2\sigma_X^2\,z_\star^3+2\sigma_X^2\, z_\star^2-2\sigma_X^2\,z_\star+\sigma_X^2 =:P\left(z_\star\right).
    \label{polynome} \end{equation}
\item The definition of $g$ in \eqref{gg} gives
     $$g(Z_n) \limas -f'(z_\star)=2\;m\; \big(z_\star \left(\mu_X-\mu_Y\right)- \mu_X \big) =-2 \,(m\sqrt{\mu_X}\sqrt{\mu_Y}),$$
    then  Proposition \ref{total-cv}  and $\alpha_n=\sqrt{1+\frac{1}{n}} -1=1/2n+O(n^{-2})$, give
    $$\frac{n}{T_{n+1}}\limas \frac{1}{m\sqrt{\mu_X}\sqrt{\mu_Y}}.$$
    By the expression of $\Gamma_{n+1}$ in \eqref{gaga}, we retrieve $\Gamma_{n+1} \limas 3/2:=\Gamma$.
\item Since the $X_n$'s are i.i.d., we claim that $X_n/\sqrt{n}\limas 0$. Indeed,  if $X_n/\sqrt{n}$ does not converge to 0, then reasoning with subsequences, we will contradict the fact that $\mu_X=\sup_n \Exp [X_n] <\infty$. The same holds for $\mu_Y$. We can then affirm that
    $$|D_{n+1}| \leq  3m\,(X_{n+1}+Y_{n+1}) \quad \mbox{and}  \quad  \Exp\big[|D_{n+1}| \;\big|\;\mathcal{F}_n\big]\leq 3\; m\;(\mu_X+ \mu_Y )  \Longrightarrow \frac{\epsilon_{n+1}}{\sqrt{n}} \limas 0 \Longrightarrow  \frac{V_{n+1}^2}{n} \limas 0.$$
    It follows that
    $$\lim_n\sum_{k=1}^n\frac{1}{n}\Exp\left[V_{k+1}^2\indi_{\left(V_{k+1}^2 \geq \varepsilon\,k\right)}\right] =\lim_n \frac{1}{n} \Exp \left[ \sum_{k=1}^{n_0}\, V_{k+1}^2\,
    \indi_{(V_{k+1}^2 \geq \varepsilon \,k)} \right]=0, \quad a.s.$$
\end{enumerate}
After these five steps, we see that Theorem \ref{th1-renlund}, then applies and we get  \eqref{tlc}. Combining last Theorem with Proposition \ref{total-cv}  and using the equality
$$\sqrt{n}\left(\frac{W_n}{n}- \frac{T_n}{n}\;z_\star\right)= \frac{T_n}{n}\;\sqrt{n} \left(Z_n-z_\star\right),$$  we obtain the last claim.
\end{proofp}
\section{Results under assumption \eqref{B}}  \label{Model 2}

In this section, we also provide a SLLN and a CLT  for the total number of ball in the urn $T_n$ under the assumption \eqref{B}. With the same notations as in the Section \ref{Model 1}, we focus on the urn model evolving under the dynamic \eqref{recusrion-opp} with $X_n=Y_n$, i.e.:
\begin{equation}\label{recurrence}\left(
    \begin{array}{c}
      W_{n+1} \\
      B_{n+1} \\
    \end{array}
  \right)\stackrel{d}{=}\left(
                          \begin{array}{c}
                            W_n \\
                            B_n \\
                          \end{array}
                        \right)+\left(
                                  \begin{array}{cc}
                                   X_{n+1}& 0 \\
                                    0 & X_{n+1}\\
                                  \end{array}
                                \right)\, \left(
                                                \begin{array}{c}
                                                  \xi_{n+1} \\
                                                  m-\xi_{n+1} \\
                                                \end{array}
                                              \right),  \quad n\geq 0.
\end{equation}
From the recursive equation \eqref{recurrence}, we have
\begin{equation}\label{tnn}
T_n=W_n+B_n=T_0+m\sum_{k=1}^{n}X_k,
\end{equation}
and the number $W_n$ of  and the proportion  $Z_n=W_n/T_n$ of white balls after $n$ draws, satisfy the recursive equations
\begin{eqnarray}
W_{n+1}&=&W_n+ X_{n+1}\,\xi_{n+1},\\
Z_{n+1}&=&Z_n\,\frac{T_n}{T_{n+1}}+\frac{X_{n+1}}{T_{n+1}}\xi_{n+1}= Z_n\, \left(1-m\,\frac{X_{n+1}}{T_{n+1}}\right)+ \frac{X_{n+1}}{T_{n+1}}\xi_{n+1}. \label{nzn}
\end{eqnarray}
\subsection{Limit theorems for  $T_n$}
The following Theorem will be useful.
\begin{teo}[Lindeberg's Theorem \cite{Lindeberg}]\label{lyp}
Let $(X_k)_k$ be a sequence of independent of square integrable random variables  defined on the same sample space. Let $s_n^2 := \sum_{k=1}^n \Vari[X_k]>0$ and assume
\begin{equation}\label{LIND}
\lim _{n\to \infty }  \frac{1}{s_n^2}\;\sum _{k=1}^n\Exp \left[(X_{k}-\Exp[X_k])^{2}\; \mathbf {1} _{\{|X_{k}-\Exp[X_k]|>\varepsilon s_{n}\}}\right]=0, \;\; \mbox{for all} \;\varepsilon >0  \quad \mbox{(Lindeberg's condition)}.
\end{equation}
Then, the CLT holds with the normalization $s_n$, i.e.
$$\frac {1}{s_n}\sum _{k=1}^n(X_k-\Exp[X_k])  \limdis N(0,1).$$
\end{teo}

Observe that the conditions of Theorem \ref{lyp} imply our condition \eqref{Model 2}. Also
observe that due to \eqref{tnn}, we
\begin{equation}\label{tv}
\Exp[T_n]= T_0+ m\sum_{k=1}^n \Exp[X_k],  \quad \Vari[T_n]=  m^2\sum_{k=1}^n \Vari[X_k],
\end{equation}
and then an immediate consequence of Theorem \ref{LIND} is:
\begin{cor}\label{total-binomial} Under Lindeberg's condition \eqref{LIND}  for the $X_n$'s of  \eqref{xy}, we have
$$\frac{T_n- \Exp[T_n]}{\sqrt{\Vari[T_n]}}\limdis {\cal N}\Big(0,\,1\Big). $$
\end{cor}
Motivated by the last CLT for the $T_n$, we provide  an SLLN in next theorem. We recall that a
function $l_0:(0,\infty)\to (0,\infty)$ is slowly varying function if it is  measurable and
\begin{equation}\label{slow}
\lim_{x\to \infty}\frac{l_0(\lambda x)}{l_0(x)}=1, \quad \forall \lambda>0.
\end{equation}
A sequence $c_n$ is regularly varying with order $\rho>0 $ if
$$c_n=n^{\rho}\, l_0(n) , \quad \mbox{and $l_0$ is slowly varying function}.$$
For the function $l_0$ we have the so-called Potter's bounds provided by  \cite[Theorem 1.5.6]{regu}: for all $ \delta>0$, there exits $A,\;N>0$ such that
\begin{equation}\label{pot}
\frac{1}{A\; n^\delta}\leq l_0(n)\leq A\; n^\delta, \quad n\geq N.
\end{equation}
Applying \cite[Proposition 1.5.8]{regu},   with $l(t)$ there equal to $l_0([t])$, we retrieve
\begin{equation}\label{regu1}
\sum_{k=1}^n c_k   \sim \frac{n\; c_n}{\rho+1}=\frac{n^{\rho+1}}{\rho+1}\, l_0(n), \quad \mbox{\it as}\; n\to \infty.
\end{equation}
\begin{teo}\label{almost sure}
Under the assumption \eqref{D}, 
\begin{equation}\label{tn}
\tilde{T}_n:=\theta_n\left(\frac{T_n}{n\; \Exp[X_n]}  -\frac{m}{\alpha+1}\right) \limas 0,
\end{equation}
for every deterministic sequence  $\theta_n$ s.t.
\begin{equation}\label{ppp}
\sum_{n\geq 1} \frac{\theta_n^2}{n^{2\alpha-\gamma+1}} \frac{l_2(n)}{ l_1(n)^2}<\infty.
\end{equation}
In particular, the latter holds for $\theta_n=1$ and $2\alpha-\gamma>0$.
\end{teo}
\begin{proofa} By \eqref{tv} and \eqref{ppp}, as $n \to\infty$,  we have
\begin{equation}\label{exab}
\Exp[T_n] = T_0+m\; \sum_{k=1}^n \Exp[X_k]  \sim \frac{m\, n}{\alpha+1}\Exp[X_n]=\frac{m}{\alpha+1} n^{\alpha+1}\, l_1(n)
\end{equation}
and
$$\Vari[T_n] =  m^2\sum_{k=1}^n \Vari[X_k]  \sim \frac{m\, n}{\gamma+1}\Vari[X_n]\sim \frac{m}{\gamma+1} n^{\gamma+1}\, l_1(n).$$
By the last asymptotic, by the condition \eqref{ppp} on $\theta_n$  and by Potter's bounds \eqref{pot} applied on $l_1$ and $l_2$, we retrieve
\begin{equation}\label{faeh1}
\frac{\Vari[T_n]}{(n\, \Exp[X_n]\,)^2} \sim  \frac{(\alpha+1)^2 }{m\, (\gamma+1)} \frac{1}{n^{2\alpha-\gamma+1}}  \frac{l_2(n)}{ l_1(n)^2} \quad \mbox{and then}\quad \sum_{n\geq 1}\frac{ \theta_n^2\;\Vari[T_n]}{(n\, \Exp[X_n]\,)^2}<\infty.
\end{equation}
Then,
$$\frac{T_n- \Exp[T_n]}{n\, \Exp[X_n]} =\frac{m}{n\, \Exp[X_n]}\sum_{k=1}^n \big(X_k-\Exp[X_k]\big)$$
and by Chebyshev's theorem, we have that for all $\varepsilon\in(0,1)$,
\begin{equation}\label{faeh2}
\Prob\left(\frac{\big|\theta_n\; (T_n- \Exp[T_n])\big|}{n\, \Exp[X_n]}\ge\varepsilon \right)\leq \left(\frac{m}{\varepsilon}\right)^2 \; \theta_n^2\;\frac{\Vari[T_n]}{(n\, \Exp[X_n])^2}
\sim   \frac{m^3}{\varepsilon^2 \;(\gamma+1)}  \; \frac{\theta_n^2\;}{n^{2\alpha-\gamma+1}} \frac{l_1(n)}{l_1(n)^2}.
\end{equation}
Finally, since $2\alpha-\gamma>0$, and due to \eqref{faeh1}, \eqref{faeh2} and by Borel-Cantelli's lemma, we conclude that
$$\lim_{n} \tilde{T}_n=\lim_{n}\theta_n \left(\frac{T_n}{n\, \Exp[X_n]} - \frac{m}{\alpha+1}\right)=\lim_{n}\theta_n\;\frac{T_n- \Exp[T_n]}{n\, \Exp[X_n]}= 0.$$
The case where $\theta_n=1$ and $2\alpha-\gamma>0$ is also justified by  Potter's bounds \eqref{pot}.
\end{proofa}
\begin{rmq}\label{r3}
\begin{enumerate}[(i)]
\item If each $X_n$ has a  Binomial distribution with parameters $\mathcal{B}(n, \, p)$,
     then $T_n$, given by \eqref{tnn}, satisfies
    $$\frac{T_n-T_0}{m} \eqdis \mathcal{B}\left(\frac{n(n+1)}{2},p\right).$$
      In this case, Lindbergs' conditions \label{LIND} are satisfied, and  also  the conditions in \eqref{D}, with $\;\alpha=\gamma=1$ and $l=p/2$ there. Hence,  Corollary \ref{total-binomial} holds.
\item By the form  \eqref{exab} of $\Exp[T_n]$, one can notice that we always have
    $$\frac{\Exp[T_n]}{m \sum_{k=1}^n\Exp[X_k]}\to l:=1+\frac{T_0}{m\;\lim_n \sum_{k=1}^n \Exp[X_k]}>0.$$
     Then, using the form the form  \eqref{exab} of $\Vari[T_n]$, we see that if  condition  \eqref{D} is replaced by
     $$ \sum_n \theta_n^2\;\frac{\sum_{k=1}^n \Vari[X_k] }{(\sum_{k=1}^n \Exp[X_k])^2}<\infty,$$
     and if we reproduce the  proof of Theorem \ref{almost sure}, then we arrive to
     $$\theta_n\;\left(\frac{T_n}{m\;\sum_{k=1}^n \Exp[X_k]}-l\right) \limas 0. $$
\end{enumerate}
\end{rmq}
\subsection{The number of  white balls under conditions \eqref{D}}

To study the asymptotic behavior of the number of white balls, we need the extra condition \eqref{D} on the sequence $(X_n)_n$  which will give more information than the almost sure convergence in \eqref{ppp}.
\begin{teo}
Under condition  \eqref{D}, the proportion $Z_n$, given by \eqref{nzn}, converges almost surely.
\end{teo}
\begin{proofa} Denoting
$$\delta_n:=\frac{X_n}{T_n}\quad \mbox{and} \quad \tilde{\xi}_n:=\xi_n-\Exp[\xi_n|\mathcal{F}_{n-1}],$$
we see that the recursive equation $Z_n$ in \eqref{nzn} takes the form
\begin{equation}\label{white}
Z_{n+1}=Z_n+\delta_{n+1}\,\tilde{\xi}_{n+1} \Longrightarrow Z_n=Z_0+\sum_{k=1}^n\delta_{k}\,\tilde{\xi}_{k}.
\end{equation}
To study the convergence of the sequence $(Z_n)_n$, we  have a problem  on the dependence of $\delta_{n+1}$ on the sigma field  $\mathcal{F}_{n+1}$, a problem that we  avoid through the following decomposition:
\begin{equation}\label{gamma}
\delta_n=\tilde{\delta}_n+\left(\delta_n-\tilde{\delta}_n\right),\quad \mbox{where} \;\; \tilde{\delta}_n:=\frac{\alpha+1}{m}\frac{X_n}{n\;\Exp[X_n]}.
\end{equation}
From the expression of $T_n$ and  $\tilde{T}_n$   given  in  \eqref{tn}, we have
$$\frac{1}{T_n}=\frac{\alpha +1}{m\;n\;\Exp(X_n)}\,\left(\frac{1}{1+\frac{\alpha +1}{m}\;\frac{\tilde{T}_n}{\theta_n}}\right),$$
then we retrieve the expressions
\begin{equation}\label{reste}
\delta_n=\frac{(\alpha +1)X_n}{m\;n\;\Exp(X_n)}\;\left(\frac{1}{1+\frac{\alpha +1}{m}\,\frac{\tilde{T}_n}{\theta_n}}\right)\quad \mbox{and} \quad \left|\delta_n-\tilde{\delta}_n\right|= \tilde{\delta}_n \;\left|\,\frac{\frac{\alpha +1}{m}\,\frac{\tilde{T}_n}{\theta_n}}{1+\frac{\alpha +1}{m}\,\frac{\tilde{T}_n}{\theta_n}}\right|.
\end{equation}
Finally, by \eqref{white} and \eqref{gamma}, we see that $Z_n$ can be decomposed into
\begin{equation}\label{decom}
Z_n=Z_0+\sum_{k=1}^n\tilde{\delta}_{k}\,\tilde{\xi}_{k}+
\sum_{k=1}^n\left(\delta_{k}-\tilde{\delta}_{k}\right) \,\tilde{\xi}_{k}=:Z_0+M_n^{(1)}+M_n^{(2)}.
\end{equation}
Now, observe  that the products  $\;\tilde{\delta}_k\,\tilde{\xi}_k, \; k=0,\;1,\;2,\ldots$, are non correlated and centered r.v.'s, and that $\tilde{\delta}_k$ is independent of $\xi_k$. We conclude that $M_n^{(1)}=\sum_{k=1}^n\tilde{\delta}_{k}\,\tilde{\xi}_{k}$ is a centred martingale which is bounded in $L^2$, i.e. $\sup_n \Exp\left[\left(M_n^{(1)}\right)^2\right]<\infty$. Actually, using $|\tilde{\xi}_{k}|\leq 2\,m$, the expression  \eqref{gamma} of $\tilde{\delta}_{k}$,  then assumption  \eqref{D} on $\Exp[X_{k}]$ and $\Vari[X_{k}]$, we obtain, as we did for obtaining  \eqref{faeh1},   that
$$\Exp\left[\left(\sum_{k\geq 1}\tilde{\delta}_{k}\,\tilde{\xi}_{k}\right)^2\right]=\Vari\left[\sum_{k\geq 1}\tilde{\delta}_{k}\,\tilde{\xi}_{k}\right]=\sum_{k\geq 1}\Vari\left[\tilde{\delta}_{k}\,\tilde{\xi}_{k}\right]\leq 4\;m^2\sum_{k\geq 1}\Vari\left[\tilde{\delta}_{k}\right]= 4\;(\alpha+1)^2   \sum_{k\geq 1} \frac{l_2(k)}{k^{2 \alpha-\gamma+2} l_1(k)^2}<+\infty.$$
Thus, $M_n^{(1)}$   converges a.s.

For the rest of this proof, we choose $\lambda\in (0,\alpha-\gamma/2)$ and $\theta_n=n^{\lambda}$.  These
choices  ensure that condition \eqref{ppp} holds and induces \eqref{tn}. Now, observing that the converging sequence  $\tilde{T}_n$ is bounded, we claim  that
\begin{equation}\label{ctrk}
 \left|\delta_k-\tilde{\delta}_k\right|   \leq \frac{\alpha+1}{m}
\;\frac{X_k}{k\;\theta_k\;\Exp[X_k]}, \quad \mbox{for $k=N,\;N+1, \ldots, \;$  and some r.v. order $N>0$.}
\end{equation}
The check the latter, use  again $|\tilde{\xi}_k|\le 2\,m$, and the fact that $|\tilde{T}_k|  \limas 0, \;|\tilde{T}_k|/\theta_k \limas 0$, hence are both smaller than 1/2 for $k$ bigger than some (random) order $N$. Then, using \eqref{reste},
we may write
$$ \left|\delta_k-\tilde{\delta}_k\right|  \;|\tilde{\xi}_k|\leq 2\;(\alpha+1) \;  \tilde{\delta}_k \;\left|\,\frac{\frac{\alpha +1}{m}\,\frac{\tilde{T}_k}{\theta_k}}{1+\frac{\alpha +1}{m}\,\frac{\tilde{T}_k}{\theta_k}}\right|\leq  \frac{4 \;(\alpha +1)^2}{m} \;\tilde{\delta}_k \; \frac{|\tilde{T}_k|}{\theta_k} \leq  \frac{2 \;(\alpha +1)^2}{m} \; \frac{\tilde{\delta}_k}{\theta_k} 
=   \frac{\alpha+1}{m} \frac{X_k}{k\;\theta_k\;\Exp[X_k]}, \quad k\geq N.$$
To conclude, use assumption \eqref{D} and then \eqref{ctrk}, write
$$\Exp\left[\sum_{k\geq 1}\frac{X_k}{k\;\theta_k\;\Exp[X_k]}\right]= \sum_{k\geq 1}\frac{1}{k^{\lambda +1}}<\infty \Longrightarrow \sum_{k\geq N}\frac{X_k}{k\;\theta_k\;\Exp[X_k]} <\infty,\; \; a.s. \Longrightarrow \sum_{k\geq 1}\left|\delta_k-\tilde{\delta}_k\right|\;|\tilde{\xi}_k| <\infty,\; \; a.s.$$
Since $M_n^{(2)}= \sum_{k=1}^n\left(\delta_k-\tilde{\delta}_k\right)\;\tilde{\xi}_k$, we deduce its $a.s.$ convergence, like the one of $M_n^{(1)}$, hence of $Z_n$.
\end{proofa}
\begin{exaample} As in Remark \ref{r3} (i),  if $X_n \eqdis B(n,\,p)$,  then all the assumptions of the last theorem are trivially satisfied.
\end{exaample}
\section{Appendix: Basic tool of stochastic approximation}\label{appendix}
The following  theorems will be used in our proofs.
\begin{teo}[Robbins-Monro, Theorem 1.4.26  pp.29, \cite{duflo}]\label{Monro}
Let $f :\mathbb{R} \to  \mathbb{R}$ be a continuous function,  and let $Z_n,\; \epsilon_n, \;\gamma_n$, be sequences of real random variable, adapted to a filtration $\mathcal{F}_n$, and  linked by the recursive equation
\begin{equation}\label{eq-monro}
Z_{n+1}=Z_{n}+\gamma_n\; \big(f(Z_n)+\epsilon_{n+1}\big), \quad n=0,\;1,\;2\ldots;
\end{equation}
Further, if the following conditions hold,
\begin{enumerate}[(a)]
\item $Z_n,\; \epsilon_n$ is square integrable, $\Exp[\epsilon_{n+1}|\mathcal{F}_n]=0$, $\gamma_n$ is nonnegative, decreases towards zero, and $\gamma_0$ is bounded by a deterministic constant;
\item $f(z^*)=0$, for some $z_\star \in \mathbb{R}$, and the scalar product $< f(z),z-z^*>$ is negative for $z\neq z_\star$;
\item $|f(z)|\leq K (1+|z|), \;\;$ for some $K>0$;
\item $\sum_n\gamma_n =\infty\;$ and $\;\sum_n\gamma_n^2\;\Exp[\epsilon_{n+1}^2|\mathcal{F}_n]< \infty,\;\; a.s.$;
\end{enumerate}
then, $Z_n \limas z^*$.
\end{teo}
\begin{teo}[Renlund \cite{renlund}]\label{th1-renlund}
Suppose $(K_n)_{n \geq 1}$ is a stochastic process adapted to a  filtration $\{\mathcal{F}_n,\; n\geq 1\}.$ Suppose it follows the recursive equation:
\begin{equation}\label{eq-renlund}K_{n+1}=\left(1-\frac{\Gamma_{n+1}}{n}\right)\,K_n+\frac{V_{n+1}}{\sqrt{n}},
\end{equation}
where, for some deterministic quantities $C,\;\sigma^2,\;\Gamma>0$, the r.v.'s $\Gamma_n,V_n \in  \mathcal{F}_n$ satisfy:
\begin{enumerate}[(a)]
\item $\Exp[V_{n+1}|\mathcal{F}_n]=o(n^{-1/2}), \; a.s.$;
\item for all $n \geq 1$, we have $\Exp[V_{n+1}^2|\mathcal{F}_n]\ \leq C, \; a.s.$
\item $\Exp[V_{n+1}^2|\mathcal{F}_n] \limas\sigma^2$;
\item $\Gamma_n\limas \Gamma$;
\item For all $\varepsilon >0$, $\sum_{k=1}^n\Exp[V_{k+1}^2\indi_{(V_{k+1}^2\geq \varepsilon k)}]\,/n\,\limas 0.$
\end{enumerate}
Then,
$$K_n\limdis\mathcal{N}\left(0,\frac{\sigma^2}{2\Gamma}\right).$$
 \end{teo}



\end{document}